\documentclass[final]{siamltex}

\usepackage{amsmath}
\usepackage{graphicx}
\usepackage{epstopdf}
\usepackage{url}

\newcommand{\half}{{\textstyle{\frac{1}{2}}}}
\newcommand{\fourth}{{\textstyle{\frac{1}{4}}}}
\newcommand{\set}[2]{\left\lbrace #1 \; : \; #2 \right\rbrace}
\newcommand{\eqnn}[1]{(\ref{#1})}
\newcommand{\pdo}[1]{\frac{\partial\phantom{#1}}{\partial{#1}}}
\newcommand{\LI}{\lambda}
\newcommand{\MI}{\mu}
\newcommand{\RI}{\rho}
\newcommand{\LGI}{\iota}

\newcommand{\bldvr}[1]{{\bf #1}}
\newcommand{\vv}{{\bf V}}
\newcommand{\dpdf}[2]{\frac{\partial{#1}}{\partial{#2}}}

\title{Optimizing the Evaluation of Finite Element Matrices}

\author{Robert C. Kirby
\thanks{Department of Computer Science, University of
Chicago, Chicago,
Illinois 60637-1581, USA }
 \and Matthew Knepley
\thanks{Division of Mathematics and Computer Science,
Argonne National Laboratory,
9700 Cass Avenue,
Argonne, Illinois 60439-4844, USA.
Work supported by the Mathematical, Information, and
Computational Sciences Division subprogram of the Office of Advanced
Scientific Computing Research,  Office of Science,
U.S. Department of Energy, under Contract
W-31-109-Eng-38
}
\and Anders Logg\thanks{Toyota Technological Institute at Chicago,
University Press Building;
1427 East 60th Street, Second Floor;
Chicago, Illinois 60637, USA;}
 \and L.~Ridgway~Scott
\thanks{The Computation Institute and
Departments of Computer Science and Mathematics, University of
Chicago, Chicago,
Illinois 60637-1581, USA.}}

\begin{document}

\maketitle

\begin{abstract}
Assembling stiffness matrices represents a significant cost in many
finite element computations. We address the question of optimizing the
evaluation of these matrices. By finding redundant computations, we
are able to significantly reduce the cost of building local stiffness
matrices for the Laplace operator and for the trilinear form for
Navier-Stokes. For the Laplace operator in two space dimensions, we
have developed a heuristic graph algorithm that searches for such
redundancies and generates code for computing the local stiffness
matrices. Up to cubics, we are able to build the stiffness matrix on
any triangle in less than one multiply-add pair per entry. Up to sixth
degree, we can do it in less than about two. Preliminary low-degree
results for Poisson and Navier-Stokes operators in three dimensions
are also promising.
\end{abstract}

\begin{keywords}
finite element, compiler, variational form
\end{keywords}

\begin{AMS}
65D05, 65N15, 65N30
\end{AMS}

\pagestyle{myheadings}
\thispagestyle{plain}
\markboth{R. C. KIRBY AND M. KNEPLEY AND A. LOGG AND L. R. SCOTT}{OPTIMIZING EVALUATION OF FINITE ELEMENT MATRICES}
\section{Introduction}
It has often been observed that the formation of the matrices arising
from finite element methods over unstructured meshes takes a
substantial amount of time and is
one of the primary disadvantages of finite elements over finite
differences.  Here, we will show that the standard algorithm for
computing finite element matrices by integration formulae is far from
optimal and present a technique that can generate algorithms with
considerably fewer operations even than well-known precomputation
techniques.

From fairly simple examples with Lagrangian finite elements, we will
present a novel optimization problem and present heuristics for the
automatic solution of this problem.  We demonstrate that the
stiffness matrix for the Laplace operator can be computed in
about one multiply-add pair per entry in two-dimensions for up to
cubics, and about two multiply-add pairs up to degree 6.  Low order examples in
three dimensions suggest similar possibilities, which we intend to
explore in the future.  More importantly,
the techniques are not limited to linear problems - in fact, we show
the potential for significant optimizations for the nonlinear term in
the Navier-Stokes equations.
These results seem to have lower flop
counts than even the best quadrature rules for simplices.

Our long-term goal is to develop a ``form compiler'' for finite
element methods.  Such a compiler will map high-level descriptions of
the variational problem and finite element spaces into low-level code
for building the algebraic systems.  Currently, the Sundance project
\cite{sundance,long} and the DOLFIN project \cite{ref:dolfin} are developing run-time C++
systems for the assembly variational forms.
Recent work in the PETSc project \cite{petsc} is leading
to compilation of variational forms into C code for building local
matrices.  Our work here complements these ideas by suggesting
compiler optimizations for such codes that would greatly enhance the
run-time performance of the matrix assembly.

Automating tedious but
essential tasks has proven remarkably successful in many areas of
scientific computing.
Automatic differentiation of numerical code has allowed
complex algorithms to be used reliably~\cite{ref:cliffetavener}.
Indeed the family of automatic differentiation tools~\cite{ref:uwe}
that
automatically produce efficient gradient, adjoint, and Hessian
algorithms for existing code have been invaluable in
enabling optimal control calculations and Newton-based nonlinear
solvers.

\section{Matrix Evaluation by Assembly}
\label{sec:evalassembl}
Finite element matrices are assembled by summing the constituent parts over each
element in the mesh.  This is facilitated through the use of a numbering scheme called the local-to-global index.
This index, $\LGI(e,\LI)$, relates the local (or element) node number,
$\LI\in{\cal L}$, on a particular
element, indexed by $e$, to its position in the global node ordering
\cite{lrsBIBfk}.  This local-to-global indexing works for Lagrangian finite elements, but requires generalizations for other families of elements.  While this generalization is required for assembly, our techniques for optimizing the computation over each element can still be applied in those situations.

We may write a finite element function $f$ in the form
\begin{equation}
  \sum_e \sum_{\LI\in{\cal L}} f_{\LGI(e,\LI)} \phi^e_\LI
\label{eqn:interpolantdef}
\end{equation}
where $f_i$ denotes the ``nodal value'' of the finite element function
at the $i$-th node in the global numbering scheme and
$\set{\phi^e_\LI}{\LI\in{\cal L}}$
denotes the set of basis functions on the element domain $T_e$.
By definition, the element basis functions, $\phi^e_\LI$,
are extended by zero outside $T_e$.
In many important cases,
we can relate all of the ``element" basis functions $\phi^e_\LI$ to a
fixed set of basis functions on a ``reference" element, ${\cal T}$, via
some mapping of ${\cal T}$ to $T_e$.
This mapping could involve changing both the ``$x$'' values and the
``$\phi$'' values in a coordinated way, as with the Piola
transform \cite{ref:piolatransfem}, or it could be one whose Jacobian
is non-constant, as with tensor-product elements or isoparametric
elements \cite{lrsBIBfk}.
But for a simple example, we assume that we have an affine mapping,
$\xi\rightarrow J\xi+x_e$, of ${\cal T}$ to $T_e$:
$$
\phi^e_\LI(x) = \phi_\LI\left(J^{-1}(x-x_e)\right).
$$
The inverse mapping, $ x\rightarrow \xi=J^{-1}(x-x_e)$ has
as its Jacobian
$$
J^{-1}_{mj} = \dpdf{\xi_m}{x_j}\,,
$$
and this is the quantity which appears in the evaluation of
the bilinear forms.
Of course, $\det J = 1/\det J^{-1}$.

The assembly algorithm utilizes the decomposition of a variational
form as a sum over ``element'' forms
$$
  a(v,w)=\sum_e a_e(v,w)
$$
where the ``element" bilinear form for Laplace's equation
is defined (and evaluated) via
\begin{equation}
\begin{split}
a_e(v,w)
 :=& \int_{T_e}\nabla v(x) \cdot \nabla w(x) \,dx   \\
  =& \int_{{\cal T}}\sum_{j=1}^d
        \pdo{x_j} v(J\xi+x_e) \\
& \pdo{x_j} w(J\xi+x_e)\det(J) \,d\xi \\
  =& \int_{{\cal T}}\sum_{j,m,m^\prime=1}^d
\dpdf{\xi_m}{x_j} \pdo{\xi_m}
\left(\sum_{\LI\in{\cal L}}v_{\LGI(e,\LI)}\phi_\LI(\xi)\right)
\times\\ &
\dpdf{\xi_{m^\prime}}{x_j} \pdo{\xi_{m^\prime}}
\left(\sum_{\MI\in{\cal L}}w_{\LGI(e,\MI)}\phi_{\MI}(\xi)
\right) \det(J)\,d\xi\\
=&  \left(\begin{matrix}
v_{\LGI(e,1)}\\ \cdot\\ \cdot\\ v_{\LGI(e,|{\cal L}|)}\\
                            \end{matrix}\right)^t
     \bldvr K^e \left(\begin{matrix}
w_{\LGI(e,1)}\\ \cdot\\ \cdot\\  w_{\LGI(e,|{\cal L}|)}\\
                            \end{matrix}\right).
\end{split}
\label{eqn:locstifalready}
\end{equation}
Here, the {\it element stiffness matrix}, $\bldvr K^e$, is given by
\begin{equation}
\begin{split}
 K^e_{\LI,\MI} :=& \sum_{j,m,m^\prime=1}^d\dpdf{\xi_m}{x_j}
                                  \dpdf{\xi_{m^\prime}}{x_j}\det(J)
    \int_{{\cal T}}\pdo{\xi_m} \phi_\LI(\xi)
                    \pdo{\xi_{m^\prime}} \phi_{\MI}(\xi) \,d\xi\\
=& \sum_{j,m,m^\prime=1}^d\dpdf{\xi_m}{x_j} \dpdf{\xi_{m^\prime}}{x_j} \det(J)
K_{\LI,\MI,m,m^\prime}\\
=&  \sum_{m,m^\prime=1}^d G^e_{m,m^\prime} K_{\LI,\MI,m,m^\prime}\\
\end{split}
\label{eqn:localeltdef}
\end{equation}
where
\begin{equation}
K_{\LI,\MI,m,m^\prime}=
    \int_{{\cal T}}\pdo{\xi_m} \phi_\LI(\xi)
                    \pdo{\xi_{m^\prime}} \phi_{\MI}(\xi) \,d\xi
\label{eqn:kaytractdef}
\end{equation}
and
\begin{equation}
G^e_{m,m^\prime}:=\det(J) \sum_{j=1}^d \dpdf{\xi_m}{x_j}\dpdf{\xi_{m^\prime}}{x_j}
\label{eqn:contractdef}
\end{equation}
for $ \LI,\MI\in{\cal L}$ and $m,m^\prime=1,\dots,d$.

The matrix associated with a bilinear form is
\begin{equation}
 A_{ij}:= a(\phi_i,\phi_j)=\sum_e a_e(\phi_i,\phi_j)
\label{eqn:defmatform}
\end{equation}
for all $i,j$, where $\phi_i$ denotes a global basis function.
We can compute this again by assembly.

First, set all the entries of $A$ to zero.
Then loop over all elements $e$ and local element numbers $\LI$
and $\MI$ and compute
\begin{equation}
\begin{split}
A_{\LGI(e,\LI),\LGI(e,\MI)} += & K^e_{\LI,\MI}
=  \sum_{m,m^\prime} G^e_{m,m^\prime} K_{\LI,\MI,m,m^\prime}\\
\end{split}
\label{eqn:compmatrxassmbl}
\end{equation}
where $G^e_{m,m^\prime}$ are defined in \eqnn{eqn:contractdef}.
One can imagine trying to optimize the computation of each
\begin{equation}
K^e_{\LI,\MI}=  \sum_{m,m^\prime} G^e_{m,m^\prime} K_{\LI,\MI,m,m^\prime}
\label{eqn:notmuchtodo}
\end{equation}
but each such term must be computed separately.
We consider this optimization in Section \ref{sec:gencomk}.

\section{Computing the Laplacian stiffness matrix with general elements}
\label{sec:gencomk}

\begin{table}
\caption{The tensor $K$ for quadratics represented as a matrix
of two by two matrices. }
\begin{center} \footnotesize
\begin{tabular}{|cc|cc|cc|cc|cc|cc|cc|cc|cc|cc|cc|cc|} \hline
  3  &  0  &  0  & -1  &  1  &  1  & -4  & -4  &  0  &  4  &  0  &  0 \\
  0  &  0  &  0  &  0  &  0  &  0  &  0  &  0  &  0  &  0  &  0  &  0 \\
\hline
  0  &  0  &  0  &  0  &  0  &  0  &  0  &  0  &  0  &  0  &  0  &  0 \\
 -1  &  0  &  0  &  3  &  1  &  1  &  0  &  0  &  4  &  0  & -4  & -4 \\
\hline
  1  &  0  &  0  &  1  &  3  &  3  & -4  &  0  &  0  &  0  &  0  & -4 \\
  1  &  0  &  0  &  1  &  3  &  3  & -4  &  0  &  0  &  0  &  0  & -4 \\
\hline
 -4  &  0  &  0  &  0  & -4  & -4  &  8  &  4  &  0  & -4  &  0  &  4 \\
 -4  &  0  &  0  &  0  &  0  &  0  &  4  &  8  & -4  & -8  &  4  &  0 \\
\hline
  0  &  0  &  0  &  4  &  0  &  0  &  0  & -4  &  8  &  4  & -8  & -4 \\
  4  &  0  &  0  &  0  &  0  &  0  & -4  & -8  &  4  &  8  & -4  &  0 \\
\hline
  0  &  0  &  0  & -4  &  0  &  0  &  0  &  4  & -8  & -4  &  8  &  4 \\
  0  &  0  &  0  & -4  & -4  & -4  &  4  &  0  & -4  &  0  &  4  &  8 \\
\hline
\end{tabular}
\end{center}
\label{tabl:geomquads}
\end{table}

The tensor $K_{\lambda,\mu,m,n}$ can be presented
as an $|{\cal L}|\times|{\cal L}|$ matrix of $d\times d$
matrices, as presented in Table \ref{tabl:geomquads} for the
case of quadratics in two dimension.
The entries of resulting matrix $K^e$ can be viewed as the dot (or Frobenius)
product of the entries of $K$ and $G^e$.
That is,
\begin{equation} \label{eqn:frobcompkay}
K^e_{\lambda,\mu}= {\bf K}_{\lambda,\mu}: G^e
\end{equation}
The key point is that certain dependencies among the entries of $K$ can be used
to significantly reduce the complexity of building each $K^e$.
For example, the four $2\times 2$ matrices in the upper-left corner
of Table \ref{tabl:geomquads} have only one nonzero entry, and six others
in $K$ are zero.  There are significant redundancies among the rest.
For example, ${\bf K}_{3,1}=-4{\bf K}_{4,1}$, so once
${\bf K}_{4,1} : G^e$ is computed,
${\bf K}_{3,1} : G^e$ may be computed by only one additional operation.

By taking advantage of these simplifications,
we see that each $K^e$ for quadratics in two dimensions
can be computed with at most 18 floating
point operations (see Section \ref{sec:computkwads})
instead of 288 floating point operations
using the straightforward definition, an improvement of a
factor of sixteen in computational complexity.

The tensor $K_{\lambda,\mu,m,n}$ for the case of linears in three dimensions
is presented in Table \ref{tabl:geomthredeepwl}.
Each $K^e$ can be computed by computing the three row sums of $G^e$,
the three column sums, and the sum of one of these sums.
We also have to negate all of the column
and row sums, leading to a total of 20 floating
point operations instead of 288 floating point operations
using the straightforward definition, an improvement of a
factor of nearly fifteen in computational complexity.
Using symmetry of $G^e$ (row sums equal column sums) we can
reduce the computation to only 10 floating
point operations, leading to a improvement of nearly 29.

\begin{table}
\caption{The tensor $K$ (multiplied by four) for piecewise linears in
three dimensions represented as a matrix of three by three matrices.}
\begin{center} \footnotesize
\begin{tabular}{|ccc|ccc|ccc|ccc|ccc|ccc|ccc|ccc|} \hline
  1  &  0  &  0  &  0  &  1  &  0  &  0  &  0  &  1  & -1  & -1  & -1 \\
  0  &  0  &  0  &  0  &  0  &  0  &  0  &  0  &  0  &  0  &  0  &  0 \\
  0  &  0  &  0  &  0  &  0  &  0  &  0  &  0  &  0  &  0  &  0  &  0 \\
\hline
  0  &  0  &  0  &  0  &  0  &  0  &  0  &  0  &  0  &  0  &  0  &  0 \\
  1  &  0  &  0  &  0  &  1  &  0  &  0  &  0  &  1  & -1  & -1  & -1 \\
  0  &  0  &  0  &  0  &  0  &  0  &  0  &  0  &  0  &  0  &  0  &  0 \\
\hline
  0  &  0  &  0  &  0  &  0  &  0  &  0  &  0  &  0  &  0  &  0  &  0 \\
  0  &  0  &  0  &  0  &  0  &  0  &  0  &  0  &  0  &  0  &  0  &  0 \\
  1  &  0  &  0  &  0  &  1  &  0  &  0  &  0  &  1  & -1  & -1  & -1 \\
\hline
 -1  &  0  &  0  &  0  & -1  &  0  &  0  &  0  & -1  &  1  &  1  &  1 \\
 -1  &  0  &  0  &  0  & -1  &  0  &  0  &  0  & -1  &  1  &  1  &  1 \\
 -1  &  0  &  0  &  0  & -1  &  0  &  0  &  0  & -1  &  1  &  1  &  1 \\
\hline
\end{tabular}
\end{center}
\label{tabl:geomthredeepwl}
\end{table}

We leave as an exercise to work out the reductions in computation
that can be done for the case of linears in two dimensions.

\subsection{Algorithms for determining reductions}
\label{sec:algreduce}

Obtaining reductions in computation can be done systematically as follows.
It may be useful to work in rational arithmetic and keep track of whether
terms are exactly zero and determine common divisors.
However, floating point representations may be sufficient in many cases.
We can start by noting which sub-matrices are zero, which ones have only
one non-zero element and so forth.
Next, we find arithmetic relationships among the sub-matrices, as follows.

Determining whether
\begin{equation} \label{eqn:firstepkayray}
{\bf K}_{\lambda,\mu}=c{\bf K}_{\lambda^\prime,\\mu^\prime}
\end{equation}
requires just simple linear algebra.
We think of these as vectors in $d^2$-dimensional space and just
compute the angle between the vectors.
If this angle is zero, then \eqnn{eqn:firstepkayray} holds.
Again, if we assume that we are working in rational arithmetic
then $c$ could be determined as a rational number.

Similarly, a third vector can be written in terms of two others
by considering its projection on the first two.
That is,
\begin{equation} \label{eqn:nexteptoodee}
{\bf K}_{\lambda,\mu}=c_1{\bf K}_{\lambda^1,\mu^1} +c_2{\bf K}_{\lambda^2,\mu^2}
\end{equation}
if and only if the projection of ${\bf K}_{\lambda,\mu}$ onto the
plane spanned by ${\bf K}_{\lambda^1,\mu^1}$ and ${\bf K}_{\lambda^2,\mu^2}$
is equal to ${\bf K}_{\lambda,\mu}$.

Higher-order relationships can be determined similarly by
linear algebra as well.
Note that any $d^2+1$ entries ${\bf K}_{\lambda,\mu}$ will be linearly
dependent, since they are in $d^2$-dimensional space.
Thus we might only expect lower-order dependences to be useful
in reducing the computational complexity.

We can then form graphs that represent the computation of $K^e$
in \eqnn{eqn:frobcompkay}.
The graphs have  $d^2$ nodes representing $G^e$ as inputs and
$|{\cal L}|\times|{\cal L}|$ nodes representing the entries of
$K^e$ as outputs.
Internal edges and nodes represent computations and temporary
storage.
The inputs to a given computation come directly from $G^e$ or
indirectly from other internal nodes in the graph.

The computation represented in \eqnn{eqn:compmatrxassmbl} would
correspond to a dense graph in which each of the input nodes is
connected directly to all $|{\cal L}|\times|{\cal L}|$ output nodes.
It will be possible to attempt to reduce
the complexity of the computation by finding sparse graphs that
represent equivalent computations.

We can generate interesting graphs by analyzing the entries of
${\bf K}_{\lambda,\mu}$ for relationships as described above.
It would appear useful to consider
\begin{itemize}
\item
entries ${\bf K}_{\lambda,\mu}$ which have only one non-zero element
\item
entries ${\bf K}_{\lambda,\mu}$ which have $2\leq k<<d^2$ non-zero elements
\item
entries ${\bf K}_{\lambda,\mu}$ which are scalar multiples of other elements
\item
entries ${\bf K}_{\lambda,\mu}$ which are linear combinations of other elements
\end{itemize}
and so forth.
For each graph representing the computation of $K^e$, we have a
precise count of the number of floating point operations, and we
can simply minimize the total number over all graphs generated by
the above heuristics.
Our examples indicate that computational simplications consisting
of an order of magnitude or more can be achieved.

We have developed a code called {\tt FErari} (for Finite Element
ReARrangemnts of Integrals) to carry out such optimizations automatically.
We describe this in detail in Section \ref{sec:ferarialg}.
This code linearizes the graph representation discussed above,
taking a particular evaluation path that is based on heuristics
chosen to approximate optimal evaluation.
We have used it to show that for classes of elements, significant
improvements in computational efficiency are available.

\subsection{Computing $K$ for quadratics}
\label{sec:computkwads}

Here we give a detailed algorithm for computing $K$ for quadratics.
Thus we have $6*K^e=$

\begin{equation} \label{eqn:kayquadalg}
\begin{pmatrix}
 3G_{11}  &  -G_{12}  &   \gamma_{11} & \gamma_{0} &   4G_{12} & 0  \\ 
 -G_{21}  &  3G_{22}  &   \gamma_{22} &     0   &  4G_{21}  & \gamma_{1} \\
\gamma_{11}&\gamma_{22}&3(\gamma_{11}+\gamma_{22})&\gamma_{0}&
                                                           0&\gamma_{1}\\
\gamma_{0} & 0 & \gamma_{0} & \gamma_2 &-\gamma_3 -8G_{22} & \gamma_3 \\ 
4G_{21} & 4G_{12} & 0 &-\gamma_3-8G_{22}&\gamma_2&-\gamma_3-8G_{11}\\ 
0 & \gamma_{1} & \gamma_{1} & \gamma_3 & -\gamma_3 -8G_{11} &\gamma_2 \\ 
\end{pmatrix}
\end{equation}
where the $G_{ij}$'s are the inputs and
the intermediate quantities $\gamma_i$ are defined and computed from
\begin{equation}\label{eqn:kquadalgint}
\begin{split}
\gamma_0=& -4\gamma_{11}, \\
\gamma_1=& -4\gamma_{22}, \\
\gamma_2=&4G_{1221}+8G_{1122}=\gamma_3+8\gamma_{12}
=8(G_{12}+\gamma_{12}), \\
\gamma_3=&4G_{1221}=4\gamma_{21} =8G_{12}
\end{split}
\end{equation}
where we use the notation $G_{ijk\ell}:= G_{ij}+G_{k\ell}$, and finally
the $\gamma_{ij}$'s are
\begin{equation} \label{eqn:dumtoilets}
\begin{pmatrix}
\gamma_{11}=G_{11}+G_{12} = G_{1112} &
\gamma_{12}=G_{11}+G_{22} = G_{1122}\\
\gamma_{21}=G_{12}+G_{21} = G_{1221} &
\gamma_{22}=G_{12}+G_{22} = G_{1222} \\
\end{pmatrix}
\end{equation}

Let us distinguish different types of operations.
The above formulas involve (a) negation, (b) multiplication of integers and
floating-point numbers, and (c) additions of floating-point numbers.
Since the order of addition is arbitrary, we may assume that the
operations (c) are commutative (although changing the order of
evaluation may change the result).
Thus we have $G_{1222}=G_{2212}$ and so forth.
The symmetry of $G$ implies that $G_{1112} = G_{1121}$ and
$G_{2122} = G_{1222}$.
The symmetry of $G$ implies that $K^e$ is also symmetric, by inspection,
as it must be from the definition.

The computation of the entries of $K^e$ procees as follows.
The computations in \eqnn{eqn:dumtoilets} are done first and require only
four (c) operations, or three (c) operations and one (b) operation
($\gamma_{21}=2G_{12}$).
Next, the $\gamma_i$'s are computed via \eqnn{eqn:kquadalgint},
requiring four (b) operations and one (c) operation. 
Finally, the matrix $K^e$ is completed, via three (a) operations,
seven (b) operations, and three (c) operations. 
This makes a total of three (a) operations,
twelve (b) operations, and three (c) operations. 
Thus only eighteen operations are required to evaluate $K^e$,
compared with $288$ operations via the formula \eqnn{eqn:notmuchtodo}.

It is clear that there may be other algorithms with the same amount
of work (or less) since there are many ways to decompose some of the
sub-matrices in terms of others.
Finding (or proving) the absolute minimum may be difficult.
Moreover, the metric for minimization should be run time, not some
arbitrary way of counting operations.
Thus the right way to utilize the ideas we are presenting may be
to identify sets of ways to evaluate finite element matrices.
These could then be tested on different systems (architectures
plus compilers) to see which is the best.
It is amusing that it takes fewer operations to compute $K^e$
than it does to write it down, so it may be that memory traffic
should be considered in an optimization algorithm that seems the
most efficient algorithm.

\section{Evaluation of general multi-linear forms}
\label{sec:trilfeval}

General multi-linear forms can appear in finite element
calculations.
As an example of a trilinear form, we consider that arising from the
first order term in the Navier-Stokes equations.  Though additional
issues arise over bilinear forms, many techniques carry over to give
efficient algorithms.

The local version of the form is defined by
\begin{equation}
\begin{split}
c_e(\bldvr u;\bldvr v,\bldvr w)
 :=& \int_{T_e}\bldvr u\cdot\nabla\bldvr v(x) \cdot \bldvr w(x) \,dx   \\
  =& \int_{T_e}\sum_{j,k=1}^d u_j(x)\pdo{x_j}v_{k}(x) w_{k}(x) \,dx   \\
  =& \int_{{\cal T}}\sum_{j,k=1}^d
        u_j(J\xi+x_e) \pdo{x_j}v_{k}(J\xi+x_e) w_{k}(J\xi+x_e)
\det(J) \,d\xi \\
  =& \int_{{\cal T}}\sum_{j,k,m=1}^d
    \left(\sum_{\LI\in{\cal L}}u_{j}^{\LGI(e,\LI)}\phi_\LI(\xi)\right)
 \dpdf{\xi_m}{x_j}
\left(\sum_{\MI\in{\cal L}}v_{k}^{\LGI(e,\MI)}\pdo{\xi_m} \phi_\MI(\xi)\right)
\times\\ & \qquad
\left(\sum_{\RI\in{\cal L}}w_{k}^{\LGI(e,\RI)}\phi_{\RI}(\xi)
\right)
\det(J) \,d\xi \\
  =& \sum_{j,k,m=1}^d \sum_{\LI,\MI,\RI\in{\cal L}}
u_{j}^{\LGI(e,\LI)} \dpdf{\xi_m}{x_j} v_{k}^{\LGI(e,\MI)} w_{k}^{\LGI(e,\RI)}
\det(J)
\times  \\ &
\int_{{\cal T}}
\phi_\LI(\xi) \pdo{\xi_m} \phi_\MI(\xi) \phi_{\RI}(\xi)
 \,d\xi \\
\end{split}
\end{equation}
Thus we have found that
\begin{equation}
\begin{split}
c_e(\bldvr u,\bldvr v,\bldvr w)
  =& \sum_{j,k=1}^d \sum_{\LI,\MI,\RI\in{\cal L}}
u_{j}^{\LGI(e,\LI)} v_{k}^{\LGI(e,\MI)} w_{k}^{\LGI(e,\RI)}
\times \\ &
 \sum_{m=1}^d \dpdf{\xi_m}{x_j} \det(J) \int_{{\cal T}}
\phi_\LI(\xi) \pdo{\xi_m} \phi_\MI(\xi) \phi_{\RI}(\xi)
\,d\xi \\
  =& \sum_{j,k=1}^d \sum_{\LI,\MI,\RI\in{\cal L}}
u_{j}^{\LGI(e,\LI)} v_{k}^{\LGI(e,\MI)} w_{k}^{\LGI(e,\RI)}
\sum_{m=1}^d \dpdf{\xi_m}{x_j} \det(J) N_{\LI,\MI,\RI,m}\\
\end{split}
\end{equation}
where
\begin{equation}
 N_{\LI,\MI,\RI,m} := \int_{{\cal T}}
\phi_\LI(\xi) \pdo{\xi_m} \phi_\MI(\xi) \phi_{\RI}(\xi)
\,d\xi
\label{eqn:navstoenndef}
\end{equation}
To summarize, we have
\begin{equation}
\begin{split}
c_e(\bldvr u,\bldvr v,\bldvr w)
  =& \sum_{j,k=1}^d \sum_{\LI,\MI,\RI\in{\cal L}}
u_{j}^{\LGI(e,\LI)} v_{k}^{\LGI(e,\MI)} w_{k}^{\LGI(e,\RI)}
N^e_{\LI,\MI,\RI,j}\\
  =& \sum_{k=1}^d \sum_{\MI,\RI\in{\cal L}}
v_{k}^{\LGI(e,\MI)} w_{k}^{\LGI(e,\RI)}
  \sum_{j=1}^d \sum_{\LI\in{\cal L}} u_{j}^{\LGI(e,\LI)} N^e_{\LI,\MI,\RI,j}\\
\end{split}
\label{eqn:sumriztrilf}
\end{equation}
where the element coefficients $N^e_{\LI,\MI,\RI,j}$ are defined by
\begin{equation}
N^e_{\LI,\MI,\RI,j} := \sum_{m=1}^d \dpdf{\xi_m}{x_j} \det(J) N_{\LI,\MI,\RI,m}.
             =: \sum_{m=1}^d \widetilde G_{mj} N_{\LI,\MI,\RI,m}
\label{eqn:nsenneedef}
\end{equation}
where $\widetilde G_{m}{j}:=\dpdf{\xi_m}{x_j} \det(J).$
Recall that $J$ is the Jacobian above,
and $J^{-1}$ is its inverse, and
$$
\left(J^{-1}\right)_{m,j} = \dpdf{\xi_m}{x_j} .
$$
Note that both $N_{\LI,\MI,\RI,(\cdot)}$ and $N_{\LI,\MI,\RI,(\cdot)}^e$
can be thought of as $d$-vectors. Moreover
$$
N_{\LI,\MI,\RI,(\cdot)}^e = \det(J) N_{\LI,\MI,\RI,(\cdot)} J^{-1}.
$$
Also note that $N_{\LI,\MI,\RI,(\cdot)} = N_{\RI,\MI,\LI,(\cdot)}$,
so that considerable storage reduction could be made if desired.

The matrix $C$ defined by $C_{ij}=c(\bldvr u,\phi_i,\phi_j)$ can be
computed using the assembly algorithm as follows.
First, note that $C$ can be written as a matrix of dimension
$|{\cal L}|\times |{\cal L}|$ with entries that are $d\times d$
diagonal blocks. In particular, let $I_d$ denote the  $d\times d$
identity matrix. Now set $C$ to zero, loop over all elements and
up-date the matrix by
\begin{equation}
\begin{split}
C_{\LGI(e,\MI),\LGI(e,\RI)}
+=& I_d \sum_{j=1}^d \sum_{\LI\in{\cal L}}
u_{j}^{\LGI(e,\LI)}  N^e_{\LI,\MI,\RI,j} \\
=& I_d \sum_{m,j=1}^d \widetilde G_{mj} \left(
 \sum_{\LI\in{\cal L}} u_{j}^{\LGI(e,\LI)} N_{\LI,\MI,\RI,m}\right) \\
=& I_d \sum_{m,\LI\in{\cal L}} \gamma_{m\lambda} N_{\LI,\MI,\RI,m}
\end{split}
\label{eqn:matrxlactd}
\end{equation}
for all $\MI$ and $\RI$, where
\begin{equation}
\gamma_{m\lambda} = \sum_{j=1}^d \widetilde G_{mj} u_{j}^{\LGI(e,\LI)}.
\label{eqn:gammtrxlctd}
\end{equation}
It thus appears that the computation of $C$ can be viewed as similar
in form to \eqnn{eqn:frobcompkay}, and similar optimization techniques
applied.
In fact, we can introduce the notation $K^{e,u}$ where
\begin{equation}
K^{e,u}_{\MI,\RI}
= \sum_{m,\LI\in{\cal L}} \gamma_{m\lambda} N_{\LI,\MI,\RI,m}
\label{eqn:newtrinotat}
\end{equation}
Then the update of $C$ is done in the obvious way with $K^{e,u}$.

\subsection{Trilinear Forms with Piecewise Linears}
\label{sec:trilplval}

For simplicity, we may consider the piecewise linear case.  Here, the
standard mixed formulations are not inf-sup stable, but the trilinear
form is still an essential part of the well-known family of stabilized
methods~\cite{hansboszepessystabilized,hughesfrancabalestra}
that do admit equal order piecewise linear
discretizations.  Moreover, our techniques could as well be used with
the nonconforming linear element~\cite{CRStokes}, which does admit an inf-sup
condition when paired with piecewise constant pressures.

In the piecewise linear case, \eqnn{eqn:navstoenndef} simplifies to
\begin{equation}
 N_{\LI,\MI,\RI,m} :=
\dpdf{\phi_\MI}{\xi_m}
\int_{{\cal T}} \phi_\LI(\xi) \phi_{\RI}(\xi) \,d\xi
\label{eqn:newnsenndef}
\end{equation}
We can think of $N_{\LI,\MI,\RI,m}$ defined from two matrices:
$N_{\LI,\MI,\RI,m}=D_{\MI,m} F_{\LI,\RI}$ where
\begin{equation}
 D_{\MI,m} :=
\dpdf{\phi_\MI}{\xi_m}=
\begin{pmatrix}
1 & 0 \\
0 & 1 \\
1 & 1
\end{pmatrix} (d=2)\quad\hbox{and}
\begin{pmatrix}
1 & 0 & 0 \\
0 & 1 & 0 \\
0 & 0 & 1 \\
1 & 1 & 1
\end{pmatrix} (d=3)
\label{eqn:deefixindef}
\end{equation}
and
\begin{equation}
 F_{\LI,\RI} :=
\int_{{\cal T}} \phi_\LI(\xi) \phi_{\RI}(\xi) \,d\xi
\label{eqn:effgramndef}
\end{equation}
The latter matrix is easy to determine.
In the piecewise linear case, we can compute integrals of
products using the quadrature rule
that is based on edge mid-points (with equal weights given
by the area of the simplex divided by the number of edges).
Thus the weights are $\omega=1/6$ for $d=2$ and $\omega=1/24$ for $d=3$.
Each of the values $\phi_\LI(\xi)$ is either $\half$ or zero,
and the products are equal to $\fourth$ or zero.
For the diagonal terms $\LI=\RI$, the product is non-zero
on $d$ edges, so $F_{\LI,\LI}=1/12$ for $d=2$ and $1/32$ for $d=3$.
If $\LI\neq\RI$, then the product $\phi_\LI(\xi) \phi_{\RI}(\xi)$
is non-zero for exactly one edge (the one connecting the corresponding
vertices), so $F_{\LI,\RI}=1/24$ for $d=2$ and $1/96$ for $d=3$.
Thus we can describe the matrices $F$ in general as having
$d$ on the diagonals, $1$ on the off-diagonals, and scaled
by $1/24$ for $d=2$ and $1/96$ for $d=3$.
Thus
\begin{equation}
\begin{split}
F=& \frac{1}{4(d+1)!}
\begin{pmatrix}
d & 1 & \cdots & 1 \\
1 & d & \cdots & 1\\
\cdot & \cdot & \cdots &\cdot \\
1 & 1 & \cdots & d\\
\end{pmatrix} \\ &
= \frac{d-1}{4(d+1)!}I_{d+1}+
\frac{1}{4(d+1)!}
\begin{pmatrix}
1 & 1 & \cdots & 1 \\
1 & 1 & \cdots & 1\\
\cdot & \cdot & \cdots &\cdot \\
1 & 1 & \cdots & 1\\
\end{pmatrix}
\end{split}
\label{eqn:whutlukslik}
\end{equation}
for $d=2$ or $3$, where $I_{d}$ denotes the $d\times d$ identity
matrix. Note that for a given $d$, the matrices in \eqnn{eqn:whutlukslik}
are $d+1\times d+1$ in dimension.

\begin{table}
\caption{The tensor $N$ (multiplied by ninety-six) for piecewise linears in
three dimensions represented as a matrix of four by three matrices.}
\begin{center} \footnotesize
\begin{tabular}{|cccc|cccc|cccc|cccc|} \hline
 3 & 1 & 1 & 1 & 0 & 0 & 0 & 0 & 0 & 0 & 0 & 0 & 3 & 1 & 1 & 1 \\
 0 & 0 & 0 & 0 & 3 & 1 & 1 & 1 & 0 & 0 & 0 & 0 & 3 & 1 & 1 & 1 \\
 0 & 0 & 0 & 0 & 0 & 0 & 0 & 0 & 3 & 1 & 1 & 1 & 3 & 1 & 1 & 1 \\
\hline
 1 & 3 & 1 & 1 & 0 & 0 & 0 & 0 & 0 & 0 & 0 & 0 & 1 & 3 & 1 & 1 \\
 0 & 0 & 0 & 0 & 1 & 3 & 1 & 1 & 0 & 0 & 0 & 0 & 1 & 3 & 1 & 1 \\
 0 & 0 & 0 & 0 & 0 & 0 & 0 & 0 & 1 & 3 & 1 & 1 & 1 & 3 & 1 & 1 \\
\hline
 1 & 1 & 3 & 1 & 0 & 0 & 0 & 0 & 0 & 0 & 0 & 0 & 1 & 1 & 3 & 1 \\
 0 & 0 & 0 & 0 & 1 & 1 & 3 & 1 & 0 & 0 & 0 & 0 & 1 & 1 & 3 & 1 \\
 0 & 0 & 0 & 0 & 0 & 0 & 0 & 0 & 1 & 1 & 3 & 1 & 1 & 1 & 3 & 1 \\
\hline
 1 & 1 & 1 & 3 & 0 & 0 & 0 & 0 & 0 & 0 & 0 & 0 & 1 & 1 & 1 & 3 \\
 0 & 0 & 0 & 0 & 1 & 1 & 1 & 3 & 0 & 0 & 0 & 0 & 1 & 1 & 1 & 3 \\
 0 & 0 & 0 & 0 & 0 & 0 & 0 & 0 & 1 & 1 & 1 & 3 & 1 & 1 & 1 & 3 \\
\hline
\end{tabular}
\end{center}
\label{tabl:tritablpl}
\end{table}

The tensor $N$ for linears in three dimensions is presented in
Table \ref{tabl:tritablpl}.
We see now a new ingredient for computing the entries of $K^{e,u}$
from the matrix $\gamma_{m,\LI}$.
Define $\gamma_m=\sum_{\LI=1}^4 \gamma_{m,\LI}$ for $m=1,2,3$,
and then $\tilde\gamma_{m,\LI}=2\gamma_{m,\LI}+\gamma_m$ for
$m=1,2,3$ and $\LI=1,2,3,4$. Then
\begin{equation} \label{eqn:kayeeeyoug}
K^{e,u}=
\begin{pmatrix}
\tilde\gamma_{11} & \tilde\gamma_{21} & \tilde\gamma_{31}&
\tilde\gamma_{11} + \tilde\gamma_{21} + \tilde\gamma_{31} \\
\tilde\gamma_{12} & \tilde\gamma_{22} & \tilde\gamma_{32}&
\tilde\gamma_{12} + \tilde\gamma_{22} + \tilde\gamma_{32} \\
\tilde\gamma_{13} & \tilde\gamma_{23} & \tilde\gamma_{33}&
\tilde\gamma_{13} + \tilde\gamma_{23} + \tilde\gamma_{33} \\
\tilde\gamma_{14} & \tilde\gamma_{24} & \tilde\gamma_{34}&
\tilde\gamma_{14} + \tilde\gamma_{24} + \tilde\gamma_{34} \\
\end{pmatrix}
\end{equation}
However, note that the $\gamma_m$'s are not computations that
would have appeared directly in the formulation of $K^{e,u}$
but are intermediary terms that we have defined for convenience
and efficiency.
This requires 39 operations, instead of 384 operations using
\eqnn{eqn:newtrinotat}.

\subsection{Algorithmic implications}
\label{sec:algimplic}

The example in Section \ref{sec:trilfeval} provides guidance
for the general case.
First of all, we see that the ``vector'' space of the evaluation
problem \eqnn{eqn:newtrinotat} can be arbitrary in size.
In the case of the trilinear form in Navier-Stokes considered
there, the dimension is the spatial dimension times the dimension
of the approximation (finite element) space.
High-order finite element approximations \cite{lrsBIBei} could
lead to very high-dimensional problems.
Thus we need to think about looking for relationship among
the ``computational vectors'' in high-dimensional spaces, e.g.,
up to several hundred in extreme cases.
The example in Section \ref{sec:trilplval} is the lowest order
case in three space dimensions, and it requires a twelve-dimensional
space for the complexity analysis.

Secondly, it will not be sufficient just to look for simple
combinations to determine optimal algorithms,
as discussed in Section \ref{sec:algreduce}.
The example in Section \ref{sec:trilplval} shows that we need
to think of this as an approximation problem.
We need to look for vectors (matrices) which closely approximate
a set of vectors that we need to compute.
The vectors $\vv_1=(1,1,1,1,0,\dots,0),
\vv_2=(0,0,0,0,1,1,1,1,0,0,0,0),
\vv_3=(0,\dots,0,1,1,1,1)$ are each edit-distance
one from four vectors we need to compute.
The quantities $\gamma_m$ represent the computations
(dot-product) with $\vv_m$.
The quantities $\tilde\gamma_{m\lambda}$ are simple perturbations
of $\gamma$ which require only two operations to evaluate.
A simple rescaling can reduce this to one operation.

Edit-distance is a useful measure to approximate the computational
complexity distance, since it provides an upper-bound on the
number of computations it takes to get from one vector to another.
Thus we need to add this type of optimization to the techniques
listed in Section \ref{sec:algreduce}.

\subsection{The FErari system}
\label{sec:ferarialg}

The algorithms discussed in Sections \ref{sec:algreduce} and
\ref{sec:algimplic} have been implemented in a prototype
system which we call FErari, for Finite Element Re-arrangement
Algorithm to Reduce Instructions.  We have used it to build optimized 
code for the local stiffness matrices for the Laplace operator using
Lagrange polynomials of up to degree six.  While
there is no perfect metric to predict efficiency of implementations due to
differences in architecture, we have taken a simple model to measure 
improvement due to FErari.
A particular algorithm could be tuned by hand, common divisors
in rational arithmetic could be sought, and more care could be taken
to order the operations to maximize register usage.
Before discussing the implementation of FErari more precisely, we point
to Table \ref{tabl:ferarisummary} to see the levels of optimization detected.
In the table, ``Entries'' refers to the number of entries in the upper triangle
of the (symmetric) matrix.  ``Base MAPs'' refers to the number of multiply
add pairs if we were not to detect dependencies at all, and ``FErari MAPs''
is the number of multiply-add pairs in the generated algorithm.  Although
we only gain about a factor of two for the higher order cases, we are
automatically generating algorithms with fewer multiply-add pairs than 
entries for linears, quadratics, and cubics.

\begin{table}
\caption{Summary of results for FErari on triangular Lagrange elements}
\begin{center} \footnotesize
\begin{tabular}{c|c||c|c}
Order & Entries & Base MAPs & FErari MAPs \\ \hline
1 & 6   & 24   & 7   \\
2 & 21  & 84   & 15  \\
3 & 55  & 220  & 45  \\
4 & 120 & 480  & 176 \\
5 & 231 & 924  & 443 \\
6 & 406 & 1624 & 867 \\
\end{tabular}
\end{center}
\label{tabl:ferarisummary}
\end{table}

FErari builds a graph of dependencies among the tensors
${\bf K}_{\lambda,\mu}$ in several stages.  We now describe each of
these stages, 
what reductions they produce in assembling $K^e$, and how much they cost
to perform.  We start by building the tensors $\left\{{\bf
  K}_{\lambda,\mu}\right\}$ for 
$1 \leq \lambda \leq \dim{P_k}$ and $\lambda \leq \mu \leq \dim{P_k}$ using FIAT
\cite{ref:FIAT}.  In discussing algorithmic complexity of our
heuristics below, we shall let $n$ be the size of this set.
Throughout, we are typically using ``greedy'' algorithms that quit
when they find a dependency.  Also, the order in which these stages
are performed matters, as if a dependency for a tensor is found at one
stage, will not mark it again later.

As we saw earlier, sometimes ${\bf K}_{\lambda,\mu} = 0$.  In this case,
$(K^e)_{\lambda,\mu} = 0$ as well for any $e$, so these entries cost nothing
to build.  Searching for uniformly zero tensors can be done
in $n$ operations.  While we found that three tensors vanish for
quadratics and cubics, none do for degrees four through six.

Next, if two tensors ${\bf K}_{\lambda,\mu}$ 
and ${\bf K}_{\lambda^\prime,\mu^\prime}$
are equal, then their dot product into $G^e$ will also be equal.  So,
we mark ${\bf K}_{\lambda^\prime,\mu^\prime}$ as depending on 
${\bf K}_{\lambda,\mu}$.
Naively, searching through the nonzero tensors is an $O(n^2)$ process,
but can be reduced to 
$O(n \log(n))$ operations by inserting the tensors into a binary tree
using lexicographic ordering of the entries of the tensors.
The fourth column of Table \ref{tabl:ferariresults}, titled ``Eq'', 
shows the number of such dependencies that FErari found.

\begin{table}
\caption{Results of FErari on Lagrange elements of degree 1 through 6}
\begin{center} \footnotesize
\begin{tabular}{c|c||c|c|c|c|c|c|c|c|c||c}
Order & Entries & Zero & Eq & Eq t & 1 Entry &Col & Ed1 & Ed2 & LC & Default & Cost \\ \hline
1 & 6 & 0 & 0 & 0 & 3 & 0 & 2 & 0 & 1 & 0 & 7\\
2 & 21 & 3 & 2 & 3 & 4 & 2 & 5 & 1 & 1 & 0 & 15\\
3 & 55 & 3 & 17 & 0 & 5 & 8 & 12 & 5 & 5 & 0 & 45\\
4 & 120 & 0 & 23 & 0 & 7 & 2 & 25 & 30 & 25 & 8 & 176\\
5 & 231 & 0 & 18 & 0 & 13 & 5 & 41 & 71 & 45 & 38 & 443\\
6 & 406 & 0 & 27 & 0 & 17 & 7 & 59 & 139 & 61 & 96 & 867\\
\end{tabular}
\end{center}
\label{tabl:ferariresults}
\end{table}

As a variation on this theme, if 
${\bf K}_{\lambda,\mu} = ({\bf K}_{\lambda^\prime,\mu^\prime})^t$, then
${\bf K}_{\lambda,\mu} : G^e =({\bf K}_{\lambda^\prime,\mu^\prime} : G^e$ since
$G^e$ is symmetric.  Hence, once $(K^e)_{\lambda,\mu}$ is computed,
$(K^e)_{\lambda^\prime,\mu^\prime}$ is free.  We may also search for these
dependencies among nonzero tensors that are not already marked as
equal to another tensor in $O(n \log(n))$ time by building a binary
tree.  In
this case, we compare by lexicographically ordering the components of
the symmetric part of each tensor (recall the symmetric part of $K$ is
$\frac{K+K^t}{2}$).  Equality of symmetric parts is necessary but not
sufficient for two matrices to be trasposes of each other.  So, if our
insertion into the binary tree reveals an entry with the same
symmetric part, we then perform an additional check to see if the two
are indeed transposes.  Unfortunately, we only found such dependencies
for quadratics, where we found three.  This is indicated in the fifth column
of Table \ref{tabl:ferariresults}, titled ``Eq t.''

So far, we have focused on finding and marking tensors that can be
dotted into $G^e$ with no work (perhaps once some other dot product
has been performed).
Now, we turn to ways of finding tensors whose dot
product with $G^e$ can be computed cheaply.  The first such way is to
find unmarked tensors ${\bf K_{\lambda,\mu}}$ with only one nonzero entry.  
This may be trivially performed in $O(n)$ time.  The sixth column of
Table \ref{tabl:ferariresults}, titled ``1 Entry,'' shows the number of 
such tensors for each polynomial degree.

If there is a constant $\alpha$ such that
${\bf K}_{\lambda^\prime,\mu^\prime} = \alpha {\bf K}_{\lambda,\mu}$, then
${\bf K}_{\lambda^\prime,\mu^\prime} : G^e = \alpha {\bf K}_{\lambda,\mu} : G^e$ and
hence may be computed in a single multiply.  Moreover, colinearity
among the remaining tensors may be found in $O(n \log(n))$ time by a
binary tree and
lexicographic ordering of an appropriate normalization.  To this end,
we divide each (nonzero) tensor by its Frobenius norm and ensure that
its first nonzero entry is positive.  If insertion into the binary
tree gives equality under this comparison, then we mark the tensors as
colinear.  The numbers of such dependencies found for each degree is
see in the ``Col'' column Table \ref{tabl:ferariresults}

Next, we seek things that are close together in edit
distance, differing only in one or two entries.  Then, the difference
between the dot products can be computed cheaply.
For each unmarked
tensor, we look for a marked tensor that is edit distance one away.
We iterate this process until no more dependencies are found, then
search for dependencies among the remainders.  We repeat this process
for things that are edit distance two apart.  In general, this process
seems to be $O(n^2)$.  Table \ref{tabl:ferariresults} has two columns, titled ``Ed1'' and ``Ed2,'' showing the
number of tensors we were able to mark in such a way.

If a tensor is a linear combination of two other tensors,
then its dot product can be computed in two multiply-add pairs.  This
is an expensive search to perform ($O(n^3)$), since we have to search
through pairs of tensors for each unmarked tensor.  However, we do
seem to find quite a few linear combinations, as seen in the ``LC'' column of 
\ref{tabl:ferariresults}.  Any remaining tensors are marked as ``Default'' in which case they are computed by the standard four multiply-add pairs.

After building the graph of dependencies, we topologically sort the
vertices.  This gives an ordering for which if vertex $u$ depends on
vertex $v$, then $v$ occurs before $u$, an essential feature to
generate code.  
We currently generate Python code for ease in debugging and integrating with the rest of our computational system, but we could just as easily generate C or Fortran.  In fact, future work holds
generating not particular code, but abstract syntax as in PETSc 3.0 \cite{petsc} as to enable code generation into multiple languages from the same
graph.
One interesting feature of the generated code is that it is completely 
unrolled -- no loops are done.  This leads to relatively large
functions, but sets up the code to a point where the compiler really
only needs to handle register allocation.  In Figure
\ref{fig:ferarilinear}, we present the generated code for computing
linears.

\begin{figure}
\caption{Generated code for computing the stiffness matrix for
linear basis functions}
\begin{verbatim}
from Numeric import zeros
G=zeros(4,"d")
def K(K,jinv):
        detinv = 1.0/(jinv[0,0]*jinv[1,1] - jinv[0,1]*jinv[1,0])
        G[0] = ( jinv[0,0]**2 + jinv[1,0]**2 ) * detinv
        G[1] = ( jinv[0,0]*jinv[0,1]+jinv[1,0]*jinv[1,1] ) * detinv
        G[2] = G[1]
        G[3] = ( jinv[0,1]**2 + jinv[1,1]**2 ) * detinv
        K[1,1] = 0.5 * G[0]
        K[1,0] = -0.5 * G[1]- K[1,1]
        K[2,1] = 0.5 * G[2]
        K[2,0] = -0.5 * G[3]- K[2,1]
        K[0,0] = -1.0 * K[1,0] + -1.0 * K[2,0]
        K[2,2] = 0.5 * G[3]
        K[0,1] = K[1,0]
        K[0,2] = K[2,0]
        K[1,2] = K[2,1]
        return K
\end{verbatim}
\label{fig:ferarilinear}
\end{figure}

\subsection{Code verification}
In general, the question of verifying a code's correctness is
difficult.  In this case, we have taken an existing Poisson solver and
replaced the function to evaluate the local stiffness matrix with
FErari-generated code.  The correct convergence rates are still
observed, and the stiffness matrices and
computed solutions agree to machine precision.  However, in the
future, we hope to 
generate optimized code for new, complicated forms where we do not
have an existing verified implementation, and the general question of
verifying such codes is beyond the scope of this present work.

\section{Computing a matrix via quadrature}
\label{sec:matquad}

The computations in equations 
(\ref{eqn:defmatform}--\ref{eqn:compmatrxassmbl})
can be computed via quadrature as
\begin{equation}
\begin{split}
A_{\LGI(e,\LI),\LGI(e,\MI)} += & 
   \sum_{\xi\in\Xi} \omega_\xi
 \nabla \phi_\LI(\xi) \cdot \left( {\bf G}^e \nabla \phi_{\MI}(\xi)
\right) \\
=&  \sum_{\xi\in\Xi} \omega_\xi
   \sum_{m,n=1}^d \phi_{\LI,m}(\xi)  G^e_{m,n} \phi_{\MI,n}(\xi) \\
=& \sum_{m,n=1}^d G^e_{m,n} 
 \sum_{\xi\in\Xi} \omega_\xi \phi_{\LI,m}(\xi) \phi_{\MI,n}(\xi) \\
=& \sum_{m,n=1}^d G^e_{m,n} K_{\LI,\MI,m,n}
\end{split}
\label{eqn:quadcompmatrxassmbl}
\end{equation}
where the coefficients $K_{\LI,\MI,m,n}$ are analogous to those
defined in \eqnn{eqn:localeltdef}, but here they are defined by quadrature:
\begin{equation}
\begin{split}
K_{\LI,\MI,m,n}
=& \sum_{\xi\in\Xi} \omega_\xi \phi_{\LI,m}(\xi) \phi_{\MI,n}(\xi) \\
\end{split}
\label{eqn:quadkaykalc}
\end{equation}
(The coefficients are exactly those of \eqnn{eqn:localeltdef} if
the quadrature is exact.)

The right strategy for computing a matrix via quadrature would 
thus appear to be to compute the coefficients $K_{\LI,\MI,m,n}$ 
first using \eqnn{eqn:quadkaykalc}, and then proceeding as before.
However, there is a different strategy associated with quadrature
when we want only to compute the {\em action} of the linear operator
associated with the matrix and not the matrix itself
(cf.~\cite{evalactpaper}).

\section{Other approaches}

We have presented one approach to optimize the computation
of finite element matrices. Other approaches have been suggested
for optimizing code for scientific computation.
The problem that we are solving can be represented as a 
sparse-matrix--times--full-matrix 
multiply, ${\cal KG}$, where the dimensions of
${\cal K}$ are, for example, $|{\cal L}|^2\times d^2$ for Laplace's
equation in $d$-dimensions using a (local) finite-element space ${\cal L}$.
For the non-linear term in Navier-Stokes, the dimensions of ${\cal K}$
become $|{\cal L}|^2\times d|{\cal L}|$.
The first dimension of ${\cal G}$ of course matches the second of ${\cal K}$,
but the second dimension of ${\cal G}$ is equal to the number of elements
in the mesh.
Thus it is worthwhile to do significant precomputation on ${\cal K}$.

At a low level, ATLAS \cite{ref:ATLAS} works with operations such 
as loop unrolling, cache blocking, etc. 
This type of optimization would not find the reductions in 
computation that FErari does, since the latter identifies 
more complex algebraic structures.
In this sense, our work is more closely aligned with that of 
FLAME \cite{ref:FLAME} which works with operations such as rank 
updates, triangular solves, etc.
The BeBOP group has also worked on related issues in sparse
matrix computation
\cite{Vuduc2004:statmodel, Vuduc2003:ata:bounds,
 Vuduc2002:smvm:bounds, Vuduc2002:sts:bounds, Im2001:sparsity:reg}.
Our work could be described as utilizing a sparsity structure in 
a matrix representation, and it can be written as multiplying a
sparse matrix times a very large set of (relatively small) vectors
(the columns of ${\cal G}$).
This is the reverse of the case considered in
the SPARSITY system \cite{Im2004:sparsity}, where multiplying a
sparse matrix times a small number of large vectors is considered.

Clearly the ideas we present here could be coupled with these
approaches to generate improved code.
The novelty of our work resides in the precomputation that is
applied to evaluate the action of ${\cal K}$.
In this way, it resembles the reorganization of computation
done for evaluation of polynomials or matrix multiply
\cite{knuthvoltwo}.

\section{Timing results}
We performed several experiments for the Poisson problem
with piecewise linear and piecewise quadratic elements.  We set up a series
of meshes using the mesh library in DOLFIN ~\cite{ref:dolfin}.  These meshes
contained between 2048 and 524288 triangles.  We timed computation of
local stiffness matrices by several techniques, their insertion into a
sparse PETSc matrix ~\cite{petsc-user-ref}, multiplying the matrix
onto a vector, and solving the linear system.  All times were observed
to be linear in the number of triangles, so we report all data as time
per million triangles. All computations were
performed on a Linux workstation with a 3 GHz Pentium 4 processor with
2GB of RAM.  All our code was compiled using {\tt gcc} with ``-O3''
optimization.

Our goal is to assess the efficacy and relevance of our proposed
optimizations.  Solver technology has been steadily improving over the
last several decades, with multigrid and other optimal strategies
being found for wider classes of problems.  When such efficient
solvers are available, the cost of assembling  the matrices
(both computing local
stiffness matrices and inserting them into the global matrix)
is much more
important.  This is especially true in geometric multigrid methods in which a
stiffness matrix can be built on each of a sequence of nested
meshes, but only a few iterations are required for convergence.
To factor out the choice of solver, we concentrate first on
the relative costs of building local stiffness matrices, inserting
them into the global matrix, and applying the matrix to a vector.

\begin{figure}
\label{fig:throughput}
\includegraphics[width=5in]{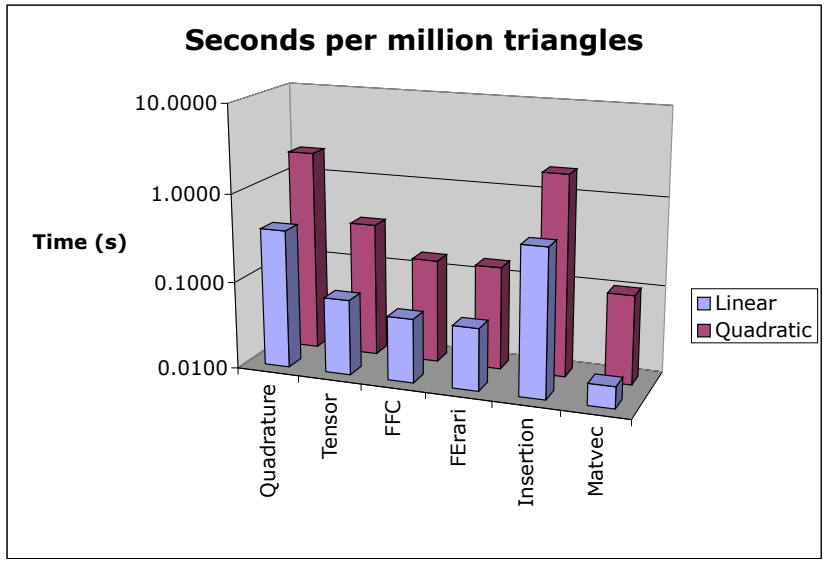}
\caption{Seconds to process one million triangles:
  local stiffness matrices, global matrix insertion, and matrix-vector product}
\end{figure}

We used four strategies for computing the local stiffness matrices -
numerical quadrature, tensor contractions (four flops per entry),
tensor contractions with zeros omitted (this code was generated by the
FEniCS Form Compiler ~\cite{ffc}), and the FErari-optimized code
translated into C.
For both linear and quadratic elements, the cost of building all of
the local stiffness matrices and inserting them into the global sparse
matrix was comparable (after storage has been preallocated).  In both
situations, computing the matrix-vector product is an order of
magnitude faster than computing the local matrices and inserting them
into the global
matrix.  These costs, all measured seconds to process one million
triangles, are given in Table~\ref{table:throughput} and plotted
Figure~\ref{fig:throughput} using log-scale for time.

\begin{table}
\label{table:throughput}

\begin{tabular}{c|cccccc}
& Quadrature & Tensor & FFC & FErari & Assemble & Matvec \\ \hline
Linear & 0.3802	& 0.0725 & 0.0535 & 0.0513 & 0.4762 & 0.0177  \\
Quadratic & 2.0000 & 0.3367 & 0.1517 & 0.1506 & 1.9342 & 0.1035
\end{tabular}
\caption{Seconds to process one million triangles:
  local stiffness matrices, global matrix insertion, and matrix-vector product}
\end{table}

We may draw several conclusions from this.  First, precomputing the
reference tensors leads to a large performance gain over numerical
quadrature.  Beyond this, additional benefits are included by omitting
the zeros, and more still by doing the FErari optimizations.  We seem
to have optimized the local computation to the point where it is
constrained by read-write to cache rather than by floating point
optimization. Second,
these optimizations reveal a new bottleneck: matrix insertion.  This
motivates studying whether we may improve the performance of insertion
into the global sparse matrix or else implementing matrix-free methods.

While the FErari optimizations are highly successful for building the
matrices, there is still quite a bit of solver overhead to contend
with.  We used GMRES (boundary conditions were imposed in a way
that broke symmetry) preconditioned by the BoomerAMG method of
HYPRE~\cite{boomeramg}.  For both linear and quadratics, the Krylov solver
converged with only three iterations.  However, the cost of building
and applying the AMG preconditioner is very large.  Using numerical
quadrature, building local stiffness matrices accounted for between
five and nine percent of the total run time (building the local to
global mapping, computing geometry tensors, computing local stiffness
matrices, sorting and inserting local matrices into the global matrix,
creating and applying the AMG preconditioner, and the rest of the
Krylov solver).  Keeping everything else constant and switching to
the FErari optimized code, building local stiffness matrices took less
than one percent of total time.  Geometric multigrid algorithms tend
to be much more efficient, and we conjecture that the cost of building
local matrices would be even more important in this case.

\section{Conclusions}

The determination of local element matrices involves a novel
problem in computational complexity.
There is a mapping from (small) geometry matrices to ``difference
stencils'' that must be computed.
We have demonstrated the potential speed-up available with
simple low-order methods.
We have suggested by examples that it may be possible to automate
this to some degree by solving abstract graph optimization problems.

\section{Acknowledgements}

We thank the FEniCS team, and Todd Dupont, Johan Hoffman, and Claes Johnson
in particular, for substantial suggestions regarding this paper.

\bibliographystyle{siam}
\bibliography{bibliography}

\end{document}